\documentclass[final,reqno]{siamltex}
\usepackage{amsmath}
\usepackage[table]{xcolor}
\usepackage{subfig}
\usepackage{graphicx}
\usepackage{amsfonts}
\usepackage{tikz}
\usepackage{placeins}
\usepackage[table]{xcolor}
\usepackage{booktabs}
\usepackage{epsfig,epsf,psfrag}
\usepackage{multirow}
\usepackage[left=1in, right=1in, top=1in, bottom=1in]{geometry}
\usepackage{algorithm}
\usepackage{algorithmic}

\title{A Kalman filter powered by $\mathcal{H}^2$-matrices for quasi-continuous data assimilation problems}
\author{Judith Y. Li \thanks{Department of Civil and Environmental Engineering, Stanford University} \and Sivaram Ambikasaran \footnotemark[4] \thanks{Courant Institute of Mathematical Sciences, New York University} \and Eric F. Darve \footnotemark[4] \thanks{Department of Mechanical Engineering, Stanford University} \and Peter K. Kitanidis \footnotemark[1] \thanks{Institute for Computational and Mathematical Engineering, Stanford University}}
\date{}

\begin{document}
\maketitle
\pagestyle{myheadings}
\markboth{Judith Y. Li, Sivaram Ambikasaran, Eric F. Darve, Peter K. Kitanidis}{A Kalman filter powered by $\mathcal{H}^2$-matrices for quasi-continuous data assimilation problems}

\begin{abstract}
Continuously tracking the movement of a fluid or a plume in the subsurface is a challenge that is often encountered in applications, such as tracking a plume of injected CO$_2$ or of a hazardous substance. Advances in monitoring techniques have made it possible to collect measurements at a high frequency while the plume moves, which has the potential advantage of providing continuous high-resolution images of fluid flow with the aid of data processing. However, the applicability of this approach is limited by the high computational cost associated with having to analyze large data sets within the time constraints imposed by real-time monitoring. Existing data assimilation methods have computational requirements that increase super-linearly with the size of the unknowns $m$. In this paper, we present the HiKF, a new Kalman filter (KF) variant powered by the hierarchical matrix approach that dramatically reduces the computational and storage cost of the standard KF from $\mathcal{O}(m^2)$ to $\mathcal{O}(m)$, while producing practically the same results. The version of HiKF that is presented here takes advantage of the so-called random walk dynamical model, which is tailored to a class of data assimilation problems in which measurements are collected quasi-continuously. The proposed method has been applied to a realistic CO$_2$ injection model and compared with the ensemble Kalman filter (EnKF). Numerical results show that HiKF can provide estimates that are more accurate than EnKF, and also demonstrate the usefulness of modeling the system dynamics as a random walk in this context.
\end{abstract}

\section{Introduction}
\label{intro}
Monitoring the progress of injected fluid fronts in the subsurface is essential to many field operations, such as storing CO$_2$ underground for greenhouse gas mitigation~\cite{benson2006monitoring,daley2007continuous},  tracking infiltration~\cite{daily1992electrical}, injecting steam for enhanced oil recovery~\cite{lazaratos1996crosswell}, and controlling hazardous substances~\cite{hubbard2001hydrogeological}. Time-lapse geophysical monitoring provides a cost-effective and noninvasive approach to image fluid flow in the volumetric region that cannot be sampled by wells. Traditional time-lapse geophysical monitoring uses a large temporal sampling interval. The vast improvements in monitoring technology allow collecting subsurface signals quasi-continuously~\cite{daley2007continuous}. Sampling the subsurface quasi-continuously, i.e., with high temporal resolution, means that we can better separate the signal (e.g., CO$_2$ leakage) from the noise, and identify important events as early as possible~\cite{arogunmati2009approach}. Therefore, a data processing method is desired to efficiently exploit the temporal redundancy resulting from fast data acquisition while also overcoming the computational challenges associated with analyzing large data sets as they arrive to obtain high-resolution fluid flow images in real time. 

One of the most appealing approaches to monitoring fluid flow in real time is to model the fluid and hydro-geophysical measurements using a state-space model (SSM) and assume the errors are Gaussian processes. Then the monitoring problems can be solved efficiently using algorithms such as Kalman filters. Kalman filter (KF) methods are powerful statistical tools for processing data as they arrive to continuously improve knowledge about a dynamical system. Recognizing that multiple solutions are consistent with the noisy data, the data conditioning process of KF is based on Bayesian inference, which combines prior knowledge and measurements to determine the posterior predictive distribution function. Unlike the Tikhonov approach~\cite{neumaier1998solving} that gives a single best solution associated with a particular choice of regularization, the solutions given by KF represent a range of possible estimates with uncertainty quantification, characterized by the maximum a posteriori (MAP) best estimate and a statistical covariance. In case of extreme events, such as CO$_2$ leakage, a proper uncertainty quantification is crucial for making informed decisions based on monitoring data.  

The original KF~\cite{Kalman:1960tn} and its nonlinear extension, the extended Kalman filter (EKF)~\cite{anderson1979optimal} do not scale well with the problem size and, thus, do not meet the high computational requirements for assimilating high frequency data. They have to store and operate on a large error covariance matrix of size $m\times m$, with $m$ in general being proportional to the number of parameters. The computational constraints limit the applications of KF to coarse-scale models, which are incapable of capturing the heterogeneity of the subsurface, and lose information on variations at fine scale, e.g., fluid flow in small faults or high-permeability channels. To address the computational challenges, one type of method reduces the effective dimension of the model state vector through applying model reduction, e.g., the approximate KF~\cite{fukumori1994approximate} and the compressed KF~\cite{guivant2001optimization}. However, those algorithms still scale quadratically with the reduced state dimension. As these model reduction techniques are generally developed for a particular application, it cannot be easily extended to work with real-time subsurface monitoring problems.   

Another type of method retains the state space, but reduces the ``working space" of KF by projecting the error covariance on a basis with reduced dimension~\cite{evensen1994sequential,TuanPham:1998hg,Pnevmatikakis:2013jh}. In other words, these algorithms resolve the computational constraints by approximating the covariance using low-rank methods, not necessarily reducing the size of the model state. For example, the singular evolutive extended Kalman filter (SEEK) approximates the error covariance matrix by a singular low rank matrix, and makes correction only in those directions associated with the singular basis. The ensemble square root filters~\cite{Tippett:2003vi}, including the ensemble Kalman filters (EnKF) and its variants~\cite{evensen1994sequential,burgers1998analysis,anderson2001ensemble,houtekamer1998data}, find the low-rank approximation of the covariance by constructing the square root of the covariance matrix using Monte Carlo methods. These methods, especially the EnKFs, have gained popularity in solving problems in hydrology and other geosciences for their computational efficiency. They do not require storing and updating the full covariance matrix, which reduces the computational cost of KF dramatically to $\mathcal{O}(m)$. 

However, replacing the full-rank covariance matrix by its low-rank approximation introduces additional errors in the estimates. It is a well-known limitation of standard Monte Carlo methods that the approximation errors decrease very slowly with the sample size, which suggests that EnKF may require a large ensemble size to obtain a solution space with small statistical error. However, the EnKF is computationally efficient only for modest sample size ($\sim 100$), which inevitably introduces large sampling errors that requires additional tuning~\cite{anderson2001ensemble} or good sampling strategies~\cite{evensen2004sampling}. The sampling bias in covariance matrices become significant when the number of observations is larger than the ensemble size~\cite{kepert2004ensemble}, which is common for assimilating geophysical data~\cite{fahimuddin2010ensemble}. Although remedies like sub-space psuedo inversion algorithm~\cite{evensen2004sampling} or sequential version of EnKF~\cite{houtekamer2001sequential} can reduce the sampling bias, the implementation of EnKF in such cases becomes less efficient.      
 
Fast linear algebra methods, such as fast fourier transform (FFT), fast multipole methods (FMM)~\cite{greengard1987fast} or hierarchical matrices approach~\cite{borm2003introduction}, provide fast and accurate alternatives to the low rank methods without reducing the rank of the covariance matrix. These methods exploit special structures in the covariance matrix that are associated with a covariance kernel, which allow them to compute the matrix-vector product very accurately, i.e. without reducing the rank of the covariance, with a near-linear scaling. For example, FFT-based inversion algorithms have been used to solve large-scale kriging problems~\cite{nowak2003efficient,Fritz:2009jn}. The $\mathcal{H}^2$ matrices, a class of hierarchical matrices~\cite{hackbusch1999sparse,hackbusch2000sparse,borm2003hierarchical} for which the FMM is applicable, can efficiently represent the covariance matrices arising out of the generalized covariance function (GCF) kernels~\cite{saibaba2012application}. Recently, the $\mathcal{H}^2$-matrices approach has been introduced to solve large-scale linear geostatistical inverse problems at $\mathcal{O} (m)$~\cite{ambikasaran2013large,ambikasaran2013fast2}. Compared to the FFT-based approach, the hierarchical matrix approach does not require uniformly-spaced regular grids~\cite{saibaba2012application}, which are inappropriate assumptions for most realistic problems. The term ``hierarchical" is a description for a certain matrix structure, which is different from the ``hierarchical filter"~\cite{anderson2007exploring}.     

In this paper, we incorporate the $\mathcal{H}^2$-matrices approach with KF to develop a computationally efficient Kalman filter (HiKF) with $\mathcal{O}(m)$ scaling. The version of HiKF that is presented here employs a random-walk forecast model, which is widely adopted in the medical-imaging literature to describe system dynamics~\cite{vauhkonen1998kalman}, and, more recently, in geophysical monitoring with rapid data acquisition~\cite{nenna2011application,quan2008stochastic}. HiKF relies on two novel ideas. First, it takes advantage of the hierarchical nature of matrices involved to accelerate the computation of dense matrix vector products. Second, the Kalman filtering equations are rewritten in a computationally efficient manner, which enables the use of these fast methods. These two new ideas reduce the computational cost of KF from $\mathcal{O} (m^2)$ to $\mathcal{O}(m)$, and allow HiKF to accurately reproduce the linear minimum mean square error (LMMSE) estimates of both the state and covariance given by the full KF. For instance, the proposed approach reduces the time of processing $10^5$ unknowns on a single core CPU from over $4$ hours to a few minutes, making it feasible for real-time data processing. Moreover, because the implementation of HiKF does not rely on any low-rank approximation of the covariance matrix, its performance will not be affected by filter inbreeding~\cite{houtekamer1998data}, a term that describes the tendency of the ensemble filters to systemically underestimate the state error covariances. 

The rest of the paper is arranged as follows. Section~\ref{section_Kalman_filter} develops the methodology for the use of Kalman filtering to solve the quasi-continuous data assimilation problems. We first introduce a state-space representation of the physical system (subsection~\ref{subsection_Linear_state}) and its solution given by the standard KF and the EnKF respectively in subsection~\ref{subsection_Kalman_filter} and~\ref{subsection_enkf}. Then the derivation of HiKF is presented in section~\ref{section:HiKF}. We first develop the idea of using the cross-covariance for error propagation (subsection~\ref{subsection:KFcross}), which is then followed by an introduction to the hierarchical matrices approach (subsection~\ref{subsection_hmatrix}) and its implementation (subsection~\ref{subsection_implementation}). In section~\ref{section_numerical_benchamrk}, we take a synthetic crosswell seismic monitoring problem to demonstrate the efficiency of HiKF, which is then compared to the standard KF and the EnKF in terms of accuracy and computational cost.
\section{Theoretical Development}
\label{section_Kalman_filter}

\subsection{Linear state space model with random walk forecast model}
\label{subsection_Linear_state}

The CO$_2$ monitoring problem can be examined in the light of estimating a hidden Markov random process. The system state of interest $\mathbf{x}_t$ is a Markov random process governed by the forecast model $f(\mathbf{x}_t \vert \mathbf{x}_{t-1})$,  where the unknown state $x_t$ is a vector of size $m$ at time step $t$. The term Markov process indicates that the value of the current state $\mathbf{x}_t$ only depends on previous state $\mathbf{x}_{t-1}$. The measurement $\mathbf{z}_t$ is a vector of size $n$, and is related to the unknown state through a measurement operator $h(\mathbf{z}_t \vert \mathbf{x}_t)$. The data assimilation problem is to recover the unobserved quantity $\mathbf{x_t}$ from a sequence of observations ${\mathbf{z_1,z_2,...,z_t}}$ that are collected at discrete time steps.

We assume that the system dynamics are governed by a linear state space model, which is given as

\begin{equation}
\mathbf{x}_{t} = F_{t} \mathbf{x}_{t-1} + \mathbf{w}_t
\label{equation_state_evolution}
\end{equation}
\begin{equation}
\mathbf{z}_{t} = H_{t}\mathbf{x}_t + \mathbf{v}_{t}
\label{equation_measurement}
\end{equation}

The state evolution equation~\eqref{equation_state_evolution} represents our knowledge about the temporal behavior of system state $\mathbf{x}_{t} \in \mathbb{R}^{m \times 1}$, where $F_{t} \in \mathbb{R}^{m \times m}$ is the state transition matrix. The input noise $\mathbf{w}_{t} \in \mathbb{R}^{m \times 1}$ is the model or input error. The measurement equation~\eqref{equation_measurement} relates the observations, $\mathbf{z}_{t} \in \mathbb{R}^{n \times 1}$ to the state vector $\mathbf{x}_{t} \in \mathbb{R}^{m \times 1}$, through the linear measurement matrix, $H_{t} \in \mathbb{R}^{n \times m}$. The vector $\mathbf{v}_{t}$ represents the measurement noise.

In many applications, the state evolution is described by a random walk model~\cite{vauhkonen1998kalman,kim20094}, in which the transition matrix $F$ is taken to be an identity matrix

\begin{equation}
\mathbf{x}_{t} = \mathbf{x}_{t-1} + \mathbf{w}_t
\label{equation_random_walk}
\end{equation}

The random-walk model is a simplification that is useful for practical applications in which data are acquired in rapid succession. Therefore, it is reasonable to assume that changes between subsequent states are small, and the incremental change can be approximated using a white-noise random process $\mathbf{w}_t$. We show later that the solution to the state-space equations can be computed very efficiently when the random-walk forecast model~\eqref{equation_random_walk} is adopted with our fast algorithm.

The noise processes $\mathbf{v}_{t}$ and $\mathbf{w}_{t}$ are specified as Gaussian processes with zero mean and known covariances given by

\begin{equation}
\mathbb{E} \left[ \mathbf{w}_t \mathbf{w}_t^T\right] = Q_t
\label{equation_state_error}
\end{equation}
\begin{equation}
\mathbb{E} \left[ \mathbf{v}_t \mathbf{v}_t^T\right] = R_t
\label{equation_measurement_error}
\end{equation}

The initial condition of the system state is assumed to be a random Gaussian variable with known mean $\mathbf{\mu}_0$ and known covariance $P_0$. That is,

\begin{equation}
\mathbb{E} \left[ \mathbf{x}_0 \right] = \mu_0
\end{equation}
\begin{equation}
\mathbb{E} \left[ (\mathbf{x}_0 - \mu_0) (\mathbf{x}_0 - \mu_0)^T\right] = P_0
\label{equation_ssm_modelnoise}
\end{equation}

As independent Gaussian processes, $\mathbf{v}_t$ and $\mathbf{w}_t$ do not depend on $\mathbf{x}_0$. The system state $\mathbf{x}_t$ and observation $\mathbf{z}_t$ are jointly Gaussian. 

\subsection{Kalman filter}
\label{subsection_Kalman_filter}
The Kalman filter~\cite{Kalman:1960tn} can be used to compute the linear minimum mean square error (LMMSE) state estimates $\mathbf{\hat{x}}_t$ compatible with the probabilistic model for the state evolution~\eqref{equation_state_evolution} and  the measurement process~\eqref{equation_measurement}. The Kalman filter is implemented in two steps: (i) Predict and (ii) Update. The prediction step estimates the state of the system at the next time step using only information from the current time step. The estimate obtained from the prediction step is an \textit{a priori} state estimate. The update step combines this prediction with the measurements obtained at the next time step to refine the estimate of the state of the system at the next time step. The estimate obtained from the update step is an \textit{a posteriori} state estimate. The \textit{a posteriori} state estimate can be thought of as a combination of the \textit{a priori} state estimate and the estimate of the system based on the new measurements.

The KF solution to the random walk state-space equations from~\eqref{equation_measurement} to~\eqref{equation_ssm_modelnoise} is given below. In all cases, we assume $R$, $Q$ and $H$ to be stationary, and hence remove the subscript $t$. Let $\hat{\mathbf{x}}_{t_2 \vert t_1}$ and $\hat{P}_{t_2 \vert t_1}$ denote the estimate of the state, $\hat{\mathbf{x}}_{t_2}$, and the covariance, $P_{t_2}$, of the system at time step $t_2$, given measurements up to time step $t_1$. 

Predict
\begin{equation}
\hat{x}_{t+1 \vert t} = \hat{x}_{t \vert t}
\end{equation}
\begin{equation} 
\hat{P}_{t+1 \vert t} = \hat{P}_{t \vert t} + Q 
\label{KF_forecast}
\end{equation}

Update
\begin{equation}
K_{t+1} = \hat{P}_{t+1 \vert t} H^T \left(H \hat{P}_{t+1 \vert t} H^T + R \right)^{-1} \\ \label{eqn_kalmangin}
\end{equation}
\begin{equation}
\hat{x}_{t+1 \vert t+1} = \hat{x}_{t+1 \vert t} + K_{t+1} \left(y_{t+1} - H \hat{x}_{t+1 \vert t} \right)
\end{equation}
\begin{equation}
\hat{P}_{t+1 \vert t+1} =\hat{P}_{t+1 \vert t} - K_{t+1} H \hat{P}_{t+1 \vert t}
\label{eqn_KF_Pa}
\end{equation}

\begin{algorithm*}
\caption{Conventional KF algorithm for a random walk forecast model ($m \gg n$)}
\label{algorithm_Kalman}
\bf{Prediction:}
\begin{center}
\rowcolors{1}{}{gray!40}
\begin{tabular}{llr}
\multicolumn{2}{c}{\bf Operation} & {\bf Cost} \\\midrule
\textit{a priori} state estimate & $\hat{x}_{t+1 \vert t} = \hat{x}_{t \vert t}$ & - \\
\textit{a priori} covariance estimate & $\hat{P}_{t+1 \vert t} = \hat{P}_{t \vert t} + Q$ & $\mathcal{O}(m^2)$\\
\end{tabular}
\end{center}

\bf{Correction:}\\
\begin{center}
\rowcolors{1}{}{gray!40}
\begin{tabular}{llr}
\multicolumn{2}{c}{\bf Operation} & {\bf Cost} \\\midrule
Kalman gain & $K_{t+1} = \hat{P}_{t+1 \vert t} H^T \left(H \hat{P}_{t+1 \vert t} H^T + R \right)^{-1}$ & $\mathcal{O}(nm^2)$\\
\textit{a posteriori} state estimate & $\hat{x}_{t+1 \vert t+1} = \hat{x}_{t+1 \vert t} + K_{t+1} \left(y_{t+1} - H \hat{x}_{t+1 \vert t} \right)$ & $\mathcal{O}(nm)$\\
\textit{a posteriori} covariance estimate & $\hat{P}_{t+1 \vert t+1} = \hat{P}_{t+1 \vert t} - K_{t+1} H \hat{P}_{t+1 \vert t}$ & $\mathcal{O}(nm^2)$
\end{tabular}
\end{center}
\end{algorithm*}
\FloatBarrier

The number of operations required for each step are summarized in Algorithm~\ref{algorithm_Kalman}. For underdetermined data assimilation problems, i.e., the number of measurements $n$ is much smaller than the number of unknowns $m$, the computational cost is measured by leading-order terms of $m$. The major computational loads of the random walk KF come from computing the Kalman gain using equation~\eqref{eqn_kalmangin} and propagating the covariance matrix using equation~\eqref{KF_forecast} and~\eqref{eqn_KF_Pa}, which all involves $\mathcal{O}(m^2)$  operations on a matrix of size $m\times m$. Therefore, in its standard form the KF is computationally prohibitive for high-dimensional modeling domains.

\subsection{Ensemble Kalman filter}
\label{subsection_enkf}
To circumvent the computational hurdle, the ensemble Kalman filter, originally proposed by~\cite{evensen1994sequential}, adopts a Monte Carlo approach to approximate the LMMSE estimates given by the KF. Begin with $N$ independent samples (ensemble) from the initial probability distribution $N(\mathbf{\bar{x}}_0,P_0)$, EnKF propagates and corrects each ensemble member $\mathbf{\tilde{x}}^j$ using the same procedures as KF. The sample mean $\mathbf{\widetilde{x}}$ and sample covariance $\widetilde{\Sigma}$ are used to approximate the mean $\mathbf{\hat{x}}$ and covariance $\Sigma$ given by the KF in section~\ref{subsection_Kalman_filter}, which are defined as

\begin{equation}
\mathbf{\hat{x}} \approx \mathbf{\tilde{x}}=\dfrac{1}{N} \sum_{j=1}^{N}\mathbf{\tilde{x}}^j
\label{equation_samplemean}
\end{equation}
\begin{equation}
\Sigma \approx \tilde{\Sigma} =\dfrac{1}{N-1} \sum_{j=1}^{N}(\mathbf{\tilde{x}}^j - \mathbf{\tilde{x}})(\mathbf{\tilde{x}}^j - \mathbf{\tilde{x}})^T = AA^T
\label{equation_samplecov1}
\end{equation}

The sample Kalman gain is computed by
\begin{equation}
\widetilde{K}_t = AA^TH^T(HAA^TH^T+R)^{-1}
\label{equation_samplegain_enkf}
\end{equation} 

The sample covariance $\widetilde{\Sigma}=AA^T$ is never constructed explicitly during the implementation, which allows the EnKF and other ensemble square root filters~\cite{Tippett:2003vi} to circumvent the computational bottleneck of Kalman filter and to be applied to high-dimensional systems. In this article we will compare our methods against the perturbation-based EnKF described in~\cite{houtekamer1998data}, the solution of which will converge to the LMMSE estimates given a large enough ensemble size $N$~\cite{burgers1998analysis}. In practice, small ensemble size limits the degrees of freedom of the sample covariance, which requires \textit{ad hoc} techniques, e.g., localization~\cite{hamill2001distance,houtekamer2001a} or covariance inflation~\cite{anderson2001ensemble} to stabilize the filter performance. Comparison of the HiKF with other versions of the EnKF, e.g., ensemble adjustment Kalman filter~\cite{anderson2001ensemble}, ensemble transform Kalman filter~\cite{bishop2001adaptive} are left for future investigations.

\section{HiKF: Kalman filter powered by $\mathcal{H}^2$ matrix approach}
\label{section:HiKF}
In this subsection, we introduce a computationally efficient Kalman filter algorithm powered by $\mathcal{H}^2$-matrices (HiKF). The computational tractability is achieved by taking advantage of a couple of key observations from Algorithm~\ref{algorithm_Kalman}. The first observation is that the initial covariance matrix can be well represented as a $\mathcal{H}^2$-matrix. This enables us to compute the matrix product $QH^T$ in $\mathcal{O}(nm)$ as opposed to standard $\mathcal{O}(nm^2)$. The second key observation is the fact that rather than store and update the full covariance matrices $\hat{P}_{t+1 \vert t}$ and $\hat{P}_{t+1 \vert t+1}$, it is enough to store and update the cross-covariance matrices, i.e., $\hat{P}_{t+1 \vert t} H^T$ and $\hat{P}_{t+1 \vert t+1} H^T$. As a result, HiKF (Algorithm~\ref{Algorithm_HiKF}) reduces the computational cost of the standard KF (Algorithm~\ref{algorithm_Kalman}) from $\mathcal{O}(m^2)$ to $\mathcal{O}(m)$. 

\begin{algorithm*}
\caption{HiKF algorithm for a random walk forecast model ($m \gg n$)}
\label{Algorithm_HiKF}
\bf{Precomputation:}
\begin{center}
\rowcolors{1}{}{gray!40}
\begin{tabular}{llr}
\multicolumn{2}{c}{\bf Operation} & {\bf Cost} \\\midrule
Cross-covariance matrix (the $\mathcal{H}^2$-matrix approach) & $C_Q = QH^T$ & $\mathcal{O}(nm)$
\end{tabular}
\end{center}

\bf{Prediction:}
\begin{center}
\rowcolors{1}{}{gray!40}
\begin{tabular}{llr}
\multicolumn{2}{c}{\bf Operation} & {\bf Cost} \\\midrule
\textit{a priori} state estimate & $\hat{s}_{t+1 \vert t} = \hat{s}_{t \vert t}$ & - \\
\textit{a priori} cross-covariance estimate & $\hat{C}_{t+1 \vert t} = \hat{C}_{t \vert t} + C_Q$ & $\mathcal{O}(nm)$
\end{tabular}
\end{center}

\bf{Correction:}
\begin{center}
\rowcolors{1}{}{gray!40}
\begin{tabular}{llr}
\multicolumn{2}{c}{\bf Operation} & {\bf Cost} \\\midrule
Kalman gain & $K_{t+1} = \hat{C}_{t+1 \vert t} \left(H \hat{C}_{t+1 \vert t} + R \right)^{-1}$ & $\mathcal{O}(n^2m)$\\
\textit{a posteriori} state estimate & $\hat{s}_{t+1 \vert t+1} = \hat{s}_{t+1 \vert t} + K_{t+1} \left(y_{t+1} - H \hat{s}_{t+1 \vert t} \right)$ & $\mathcal{O}(nm)$\\
\textit{a posteriori} cross-covariance estimate & $\hat{C}_{t+1 \vert t+1} = \hat{C}_{t+1 \vert t} - K_{t+1} H \hat{C}_{t+1 \vert t}$ & $\mathcal{O}(n^2m)$\\
\textit{a posteriori} variance estimate & $\delta^2_{t+1 \vert t+1} = \delta^2_{t+1 \vert t} - \sum\limits_{j=1}^{n}(K_{t+1})_{ij} (\hat{C}_{t+1 \vert t})_{ij}$ & $\mathcal{O}(nm)$
\end{tabular}
\end{center}
\end{algorithm*}
\FloatBarrier

\subsection{Uncertainty propagation using cross-covariances}
\label{subsection:KFcross} 
In essence, KF requires storing and updating the covariance matrix $P$ between adjacent steps, which results in a quadratic computational cost. Computing Kalman gain is also expensive as it requires forming the cross-covariance matrix $C = PH^T$ from the covariance matrix $P$. However, if we assume the prior cross-covariance matrix $\hat{C}_{t+1|t}$ has already been computed, then Kalman gain can be computed at the cost of $\mathcal{O} (n^2m)$ according to equation~\eqref{eqn_kalmangin}, 
\begin{equation}
K_{t+1} = \hat{C}_{t+1 \vert t} \left(H \hat{C}_{t+1 \vert t} + R \right)^{-1} 
\label{eqn_kalmangin_new}
\end{equation}

where the prior cross-covariance can be computed from the cross-covariance from last step and the product of $\hat{C}_Q=QH^T$,
 \begin{equation}
\hat{C}_{t+1 \vert t} = \hat{C}_{t \vert t} + \hat{C}_Q
\label{eqn_KF_forecast_new}
\end{equation}

Equation~\eqref{eqn_KF_forecast_new} is based on the assumption of a stationary observation operator $H$ and the random walk dynamical model~\eqref{equation_random_walk}. Given a fixed observation operator $H$ and the computed Kalman gain, the cross-covariance $C$ can be updated at $\mathcal{O}(n^2m)$ according to variations of equation~\eqref{KF_forecast} and~\eqref{eqn_KF_Pa}, that is, 
\begin{equation}
\hat{C}_{t+1 \vert t+1} =\hat{C}_{t+1 \vert t} - K_{t+1} (H \hat{C}_{t+1 \vert t})
\label{eqn_KF_pred_new}
\end{equation}

The posterior variance is given by
\begin{equation}
\delta^2_{t+1 \vert t+1} = \delta^2_{t+1 \vert t} - \sum\limits_{j=1}^{n}(K_{t+1})_{ij} (\hat{C}_{t+1 \vert t})_{ij} \end{equation}

Therefore, we have derived an equivalent formulation of the KF by substituting equation~\eqref{KF_forecast} to~\eqref{eqn_KF_Pa} with equation~\eqref{eqn_kalmangin_new} to~\eqref{eqn_KF_pred_new}. With the assumptions that $QH^T$ and $\hat{C}_{0 \vert 0}$ have already been calculated, the new formulation obviates the need to evaluate and store the full covariance matrix $P$. Instead, the cross-covariance $C \in \mathcal{R}^{m\times n}$ from last time step is stored and updated at each data assimilation step in order to predict the state estimate at current time step. Therefore, the cost for each measurement update is $\mathcal{O}(m)$ instead of $\mathcal{O}(m^2)$.

To initialize the new KF algorithm, one needs the initial state vector $\hat{x}_{0|0}$, its error cross-covariance $\hat{C}_{0 \vert 0}= P_{0|0}H^T$ and the model error cross-covariance $C_Q = QH^T$. Since after a few iterations the estimation process is usually independent from the choice of $P_{0|0}$, the initial covariance $P_{0|0}$ is assumed to be $\alpha I_m$, where $I_m$ is a $m\times m$ identity matrix, and $\alpha$ can be chosen to be an arbitrary large value if no information about the value of $P_{0|0}$ is available. Hence the cross-covariance of the initial state, i.e., $\hat{C}_{0 \vert 0}= P_{0|0}H^T$ can be initialized to be $\alpha H^T$. The model error covariance $Q$ is a $m\times m $ matrix, with each element $Q_{ij}$ defined using the generalized covariance function (GCF) $K(r_{ij})$, where $r_{ij}$ is the spatial distance between $x_i$ and $x_j$, to account for the spatial variability~\cite{pk1997introduction}. We will show in next section that this kernel representation gives $Q$ an underlying $Hierarchical$ structure. Therefore, despite that $Q$ is a dense matrix, the matrix product $C_Q = QH^T$ can be computed at the cost of $\mathcal{O}(m)$ as opposed to $\mathcal{O}(m^2)$.     

\subsection{Hierarchical Matrices}
\label{subsection_hmatrix}
The $\mathcal{H}$-matrix approach provides a data-sparse representation of dense matrices arising in many engineering applications, e.g., the boundary element method (BEM) for discretizing integral equations~\cite{bebendorf2000approximation} and dense covariance matrices~\cite{ambikasaran2013large}. The data-sparse representation relies on the fact that these matrices can be recursively sub-divided based on a tree structure, and most matrix sub-blocks at different levels in the tree can be well-approximated using a low-rank block. There are many different hierarchical matrices depending on the tree structure and algorithms to obtain the low-rank blocks. The main advantage of the $\mathcal{H}$-matrix approach is that the storage cost and the computational cost of matrix algebra operations of a $N \times N$ $\mathcal{H}$-matrix scales as $\mathcal{O}(N \log^a N)$ instead of $\mathcal{O}(N^2)$ or $\mathcal{O}(N^3)$~\cite{hackbusch1999sparse, hackbusch2000sparse, bebendorf2008hierarchical, borm2003introduction, ambikasaran2013fast,martinsson2005fast,chandrasekaran2006fast,chandrasekaran2006fast1}.

In this article, we will restrict our attention to a class of hierarchical matrices known as $\mathcal{H}^2$ matrices, which efficiently represent the covariance matrices arising out of the GCF kernels~\cite{ambikasaran2013large}. $\mathcal{H}^2$-matrices are the subset of the $\mathcal{H}$ matrices for which the fast multipole method (FMM) is applicable. One characteristic of an error covariance Q with $\mathcal{H}^2$ structure is that the blocks corresponding to a distant cluster of points has low rank even though the full matrix is dense. Observe that if $H^T$ is a matrix and we wish to evaluate $QH^T$, then we have
\begin{equation}
QH^T = [Qh^{(1)}\,Qh^{(2)}\, ...\, Qh^{(n)} ]
\end{equation}
Hence, the FMM~\cite{greengard1987fast, fong2009black} can be employed to compute and store the matrix-vector products of the summation form
\begin{equation}
Qh_i^{(k)} = \sum_{j=1}^m K(x_i,x_j)h_j^{(k)}
\end{equation}
for $i \in \{1,2,...,m\}$ in $\mathcal{O}(m)$ operations as opposed to $\mathcal{O}(m^2)$ with a controllable error $\epsilon$.

A key part of the FMM is the use of fast low-rank factorization algorithms to represent the interactions $K(x_i,x_j)$ between well-separated clusters of points. There are several low-rank factorization techniques available. For instance, if the entries of the kernel matrix $K$ arise from a Green's function or some analytic function (in some suitable region), the low-rank approximation can be computed using the multipole expansion, as in the traditional fast multipole method~\cite{greengard1987fast}, or using other analytic techniques such as Taylor series expansions~\cite{borm2003hierarchical} or Chebyshev interpolation~\cite{fong2009black}. However, if the entries in the matrix $K$ are only known algebraically, then fast algebraic techniques, including adaptive cross approximation (ACA)~\cite{rjasanow2002adaptive} and pseudo-skeletal approximations~\cite{goreinov1997theory}, can be used to obtain low-rank factorizations in almost linear complexity. In this work, we rely on Chebyshev interpolation~\cite{fong2009black,ambikasaran2013fast} to construct fast low-rank factorizations. Theorems are available~\cite{fong2009black} to determine the rank necessary to achieve accuracy $\epsilon$.

\subsection{Implementation}
\label{subsection_implementation}
A summary of the HiKF algorithm can be found in Algorithm~\ref{Algorithm_HiKF}. The overall cost of the HiKF is $\mathcal{O}(n^2m)$ as opposed to $\mathcal{O}(nm^2)$. Recall that for most problems $n \ll m$. Therefore, the solutions can be obtained in $\mathcal{O}(m)$ time. C++ Software packages related to the fast methods used in this article, including FLIPACK and BBFMM (black box FMM), can be operated as a block box and are available online at https://github.com/sivaramambikasaran/. A MATLAB version of BBFMM is available at https://github.com/judithyueli/mexBBFMM2D, which allows the readers to implement HiKF in MATLAB. We refer readers to~\cite{ambikasaran2013large} and~\cite{fong2009black} for details of the black box FMM algorithm. It is also worth noting that~\cite{ambikasaran2013large} also exploits the sparsity of the measurement operator $H$ in reducing the run time. The black box FMM is compatible with a wide range of kernels including GCFs~\cite{kitanidis1993generalized,kitanidis1983statistical,kitanidis1999generalized,starks1982estimation} . Some of the commonly used kernels are listed below,
\begin{enumerate}
\item Gaussian kernels
\begin{equation}
Q(r) = \sigma ^2 \exp (-\dfrac{r^2}{L^2})
\end{equation}
\item Exponential kernels
\begin{equation}
Q(r) = \exp(-\dfrac{r}{L})
\end{equation}
\item Logrithm kernels
\begin{equation}
Q(r) = A\log (r),A<0
\end{equation}
\end{enumerate}

A chart that demonstrates the scalability and the accuracy of BBFMM for the Gaussian kernel is shown in Table~\ref{table_bbfmm}. As shown from the table, usually a highly accurate approximation ($\approx 10^{-9}$) can be obtained using around 5 Chebyshev nodes. The accuracy can be improved by increasing the number of Chebyshev nodes with only a slight increase in computational time. BBFMM has linear scaling and it takes only two minutes to compute the matrix-vector product that contains one-million entries.

\begin{center}
\begin{table}
\center
\caption{Comparison of BBFMM and the standard approach on matrix-vector multiplication using a 2.40 GHz single-core CPU}
\small{
\rowcolors{1}{}{gray!40}
\begin{tabular}{c|c|c|c|c}
\hline
\bf{State dimension} & \multicolumn{2}{c}{\bf{Time (seconds)}} \vline & \multicolumn{1}{c}{\bf{Relative accuracy}} \vline & \multicolumn{1}{c}{\bf{Number of Chebyshev nodes}}\\
\hline
& Standard & BBFMM &   &  \\
\hline  
$10,000$ & $7.25$ & $0.61$ & $4.08e-9$ & 5  \\
 &  & $0.96$  & $5.84e-10$& 6  \\
 &  & $1.07$ & $3.28e-11$ & 7  \\
$100,000$&$1627$&$11.54$&$3.27e-11$&$7$\\
$1,000,000$& $-$ &$132.6$&$-$&$7$\\
\end{tabular}
}
\label{table_bbfmm}
\end{table}
\end{center}
\FloatBarrier

\section{Numerical Benchmark}
\label{section_numerical_benchamrk}
In this section, we illustrate the performance of HiKF by comparing it against the conventional KF and EnKF. The data assimilation methods  are implemented to continuously track a CO$_2$ plume in the subsurface from seismic travel times, resulting in a dynamical system described by a linear state space model with random walk forecast model.
 
\subsection{Synthetic TOUGH2 CO$_2$ Monitoring Experiment}
\begin{figure}
\centering\includegraphics{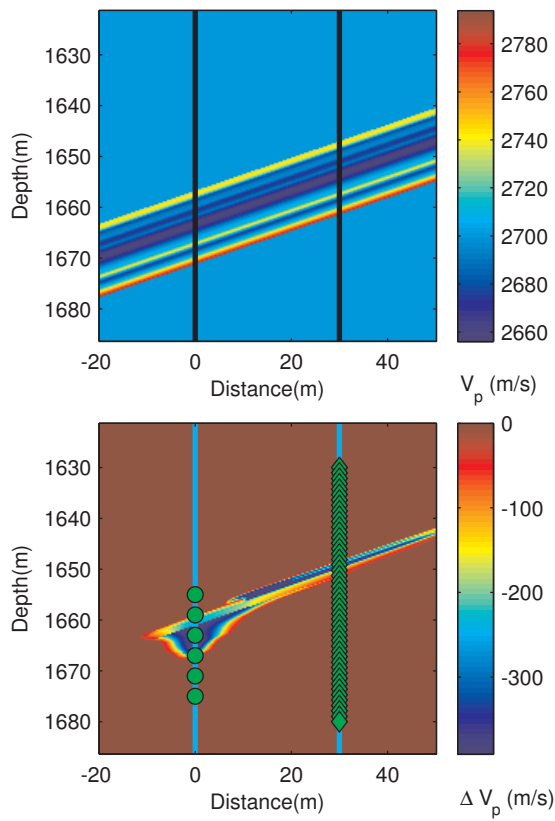}
\caption{The base velocity model (top) before CO$_2$ injection and the simulated CO$_2$-induced velocity reduction (bottom) 5 days after injection at FRIO II site. Data provided courtesy of Jonathan B. Ajo-Franklin, Thomas M. Daley and Christine Doughty from Earth Science Division, Lawrence Berkeley National Laboratory, Berkeley, CA, United States.}
\label{figure_BaseSurvey}
\end{figure}

A detailed reservoir flow model has been developed for the Frio-II brine pilot CO$_2$ injection experiment~\cite{daley2007continuous} using a state-of-the-art model reservoir simulator TOUGH2~\cite{pruess1991tough2}. The Frio-II pilot injected 380 tons of CO$_2$ into a $17$m-thick deep brine aquifer at the depth of $1657$m. Based on well logs and core measurements, the brine aquifer is assumed to have a dip of $18$ degree and a layered permeability and porosity distribution. The flow simulation predicted the spatial distributions of CO$_2$ saturation and pressure for $5$ days. The seismic wave velocity model is built from the flow model using a patchy petrophysical relationship~\cite{White:2001ul}. The baseline and the time-lapse changes of the compressional seismic wave velocities are shown in Figure~\ref{figure_BaseSurvey}. CO$_2$ can be monitored seismically by mapping the time-varying CO$_2$-induced velocity reductions from measurements of traveltime delays~\cite{wang1998seismic}. We conduct synthetic crosswell seismic surveys every $3$ hours, with $6$ sources deployed at the injection well and $48$ receivers deployed at the observation well $30$ meters apart. The acquisition geometry follows~\cite{ajo2009optimal} and remains fixed during the monitoring experiment.  

\FloatBarrier

The travel time measurements are simulated by integrating the seismic slowness (reciprocal of velocities) along the raypaths on which seismic waves propagate. In reality, the ray path will change in response to CO$_2$ plume migration. For simplicity, we assume that the raypath follows a straight line connecting the source-receiver pair, which is the high-frequency limiting case of the wave equation. By assuming the raypath $l$ has not been significantly altered due to CO$_2$ injection, the traveltime $\mathbf{y}_t$ can be expressed as a linear function of the unknown slowness $\mathbf{s}_t$, which varies with location $\mathbf{r}$. That is,

\begin{equation}
\mathbf{y}_t = \int _l \mathbf{s}(\mathbf{r})dl = H\mathbf{s}_t
\label{equation_traveltimeintegral}
\end{equation}

The time-invariant linear observation operator $H$ embeds the information about the area sampled by the raypath, with each element $H_{ij}$ representing the length of the $i^{th}$ raypath within the $j^{th}$ cell. The CO$_2$-induced low velocity zone is imaged from traveltime delays $\Delta \mathbf{y_t}$ relative to the baseline traveltime, which can be obtained by subtracting the baseline traveltime from the measurement equation~\eqref{equation_traveltimeintegral}. That is,

\begin{equation}
\Delta \mathbf{y}_t = H\Delta \mathbf{s}_t
\label{equation_difftomo} 
\end{equation}

The variable of interest $\Delta \mathbf{s}_t$ is the perturbation of the background slowness at time step $t$. The differential tomography approach applies spatial and temporal regularizations directly on the slowness perturbations instead of slowness itself. The observation $\mathbf{z}_t$ in equation~\eqref{equation_measurement} is simulated as traveltime delay $\Delta \mathbf{y}_t$ contaminated with white noise $\mathbf{v}_t$, resulting in a 65dB signal-to-noise ratio (SNR) defined as

\begin{equation}
SNR = 10\log_{10} \dfrac{\Vert {\Delta \mathbf{y}} \Vert^2_2}{\Vert \mathbf{\sigma} \Vert^2_2}
\label{equation_SNR}
\end{equation}

where $\Delta \mathbf{y}$ is the  measurement signal, and $\mathbf{\sigma ^2}$ is the measurement variance. The measurement noise $\mathbf{v}$ are realizations from $N(0,\sigma ^2I)$.

\subsection{State-space model parameter selection}

The state-space representation of the dynamic traveltime tomography problem is given in equation \eqref{equation_measurement} to \eqref{equation_ssm_modelnoise}. The major error sources in the mathematical model come from using the random-walk model \eqref{equation_random_walk} to approximate the true complex system dynamics, as well as using noisy data sets. If the data assimilation method is based on Monte Carlo, e.g. EnKF, then sampling errors become another error source. The parameterization of the noise structures incorporates our prior knowledge of the dynamical process, which is given by

\begin{equation}
R = \sigma^2 I
\label{equation_R}
\end{equation} 
\begin{equation}
Q_{ij} = \theta \exp(-\dfrac{||\mathbf{x}_i-\mathbf{x}_j|| ^p}{l^p})
\label{equation_Q}
\end{equation}

The observation noise $\mathbf{v}_t\sim N(0,R)$ is assumed to be uncorrelated, where $I\in \mathbb{R}^{n\times n}$ is an identity matrix, and $\sigma ^2$ is the observation variance. The noise process $\mathbf{w}_t \sim N(0,Q)$ accounts for the lag errors due to using a random-walk forecast model. As shown in equation~\eqref{equation_Q}, the model error covariance $Q$ is parameterized as a kernel function, which possesses a hierarchical structure described in section~\ref{subsection_hmatrix}. The correlation between two points decays as their spatial separation increases. Small power $p$ indicates sharp changes and a large $l$ indicates long-range correlation. The lag errors can be estimated approximately by selecting reasonable values for $p$ and $l$. The regularization parameters $\theta$ and $\sigma$ control the relative influence of the lag errors and the observational errors. If the slowness changes are significant, then the lag errors can be assumed to prevail, therefore more weight should be put on measurements. Although $Q$ is a large full-rank matrix, it can be efficiently represented in a data-sparse manner as a hierarchical matrix, more specifically, as a $\mathcal{H}^2$ matrix~\cite{ambikasaran2013large,saibaba2012application}. Therefore, the product $QH^T$ can be computed at $\mathcal{O}(mn)$, where $m \gg n$. 
    
As we assume no CO$_2$ is present before the injection, the initial guess for the slowness perturbation $\mathbf{x}_0$ is a zero vector, with zero initial error covariance $P_0$. The effects of the choice of $P_0$ die out relatively fast as more data get assimilated~\cite{vauhkonen1998kalman}.

\subsection{Dynamic inversion results}
    
The slowness changes induced by CO$_2$ injection have been estimated from simulated traveltimes acquired at 41 time frames using KF, HiKF and EnKF described in previous sections. All prediction filters adopt the same random-walk state space equations from \eqref{equation_measurement} to \eqref{equation_ssm_modelnoise} to describe the dynamical system, with the same noise structures and regularization parameters. Thus, the difference in the estimates results solely from the choice of filter. Methods are compared in terms of accuracy, uncertainty quantification, and computational cost. CO$_2$-induced slowness changes are estimated at three different resolutions, with $59\times 55$, $117\times 109$, and $234\times 217$ pixels, respectively. The accuracy of the inversion results is measured in terms of relative estimation error $\mathbf{e}_t$ defined as

\begin{equation}
\mathbf{e}_t = \dfrac{\left \Vert \mathbf{x}_{\text{est}}-\mathbf{x}_{\text{ref}} \right \Vert _2}{\left \Vert \mathbf{x}_{\text{ref}} \right \Vert _2}
\label{equation_RelativeEstErr}
\end{equation}

\begin{figure}
\centering\includegraphics[scale=0.8]{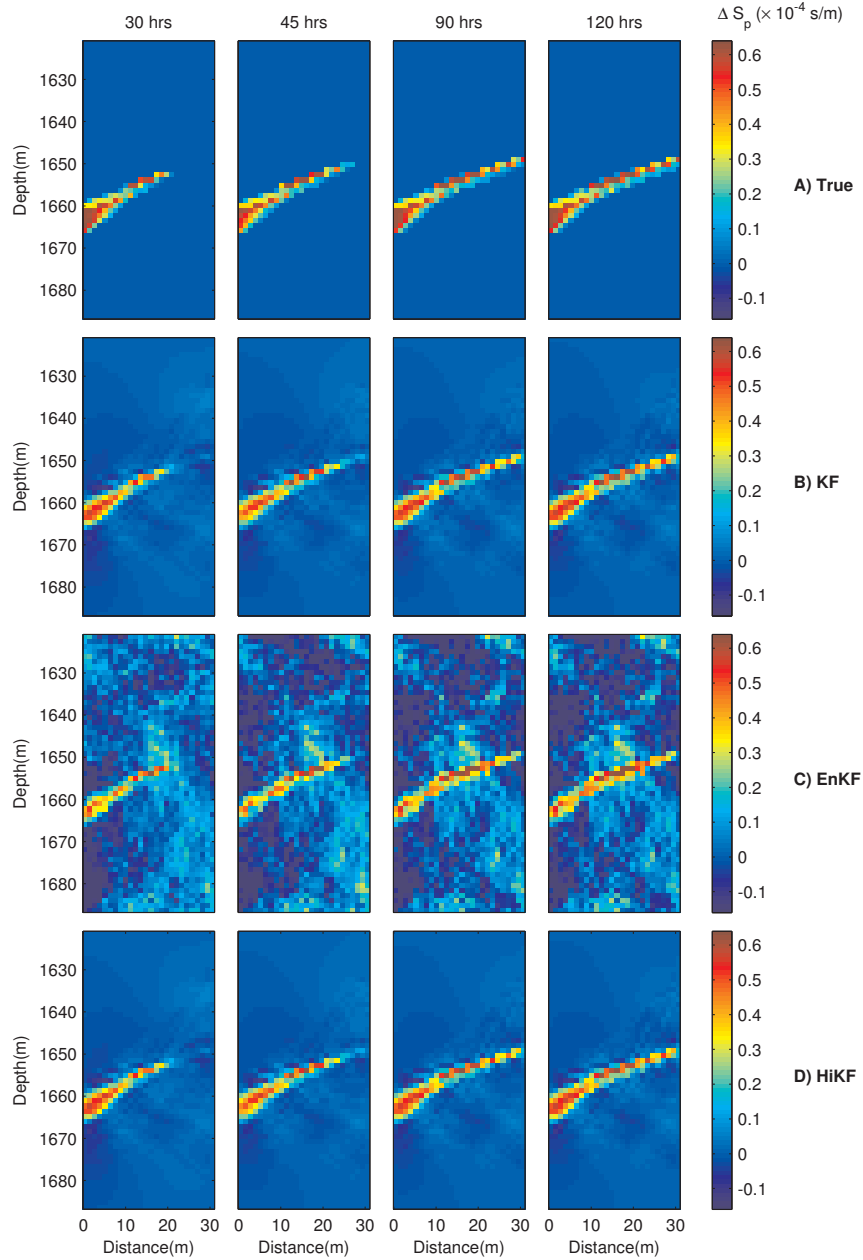}
\caption{True and estimated CO$_2$-induced changes in compressional-wave slowness ($\times 10^{-4}m/s$) between two wells at low resolution. EnKF results are obtained with 600 samples. Each state vector $\mathbf{x}$ represents a $59\times 55$ image, each observation vector $\mathbf{z}$ contains 288 simulated tomographic travel time measurements with 65-dB SNR at each time instance. SNR is defined in \eqref{equation_SNR}. Data is acquired at 41 time instances, lasting for 5 days.}
\label{figure_EnKFest_59x55}
\end{figure}

\begin{figure}
\centering\includegraphics{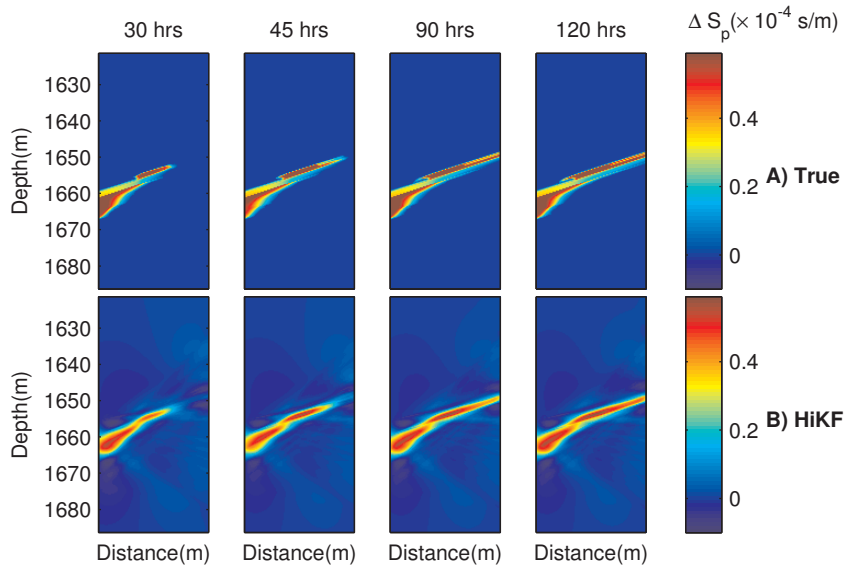}
\caption{True and HiKF estimates of the compressional-wave slowness ($\times 10^{-4}s/m$) at high resolution. Each state vector $\mathbf{x}$ represents a $234\times 217$ image.}
\label{figure_EnKFest_234x217}
\end{figure}

\begin{figure}
\centering\includegraphics{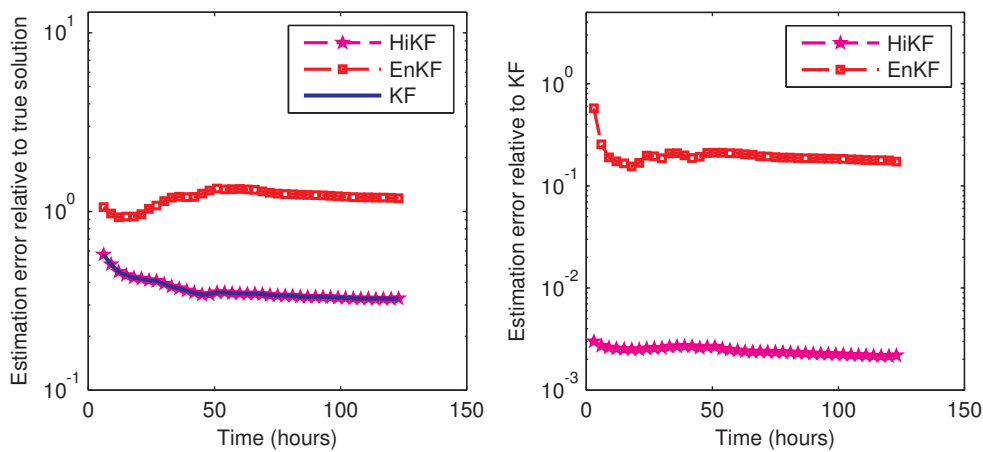}
\caption{Estimation error relative to the true solution (left) and to the LMMSE solution by KF (right). Definition is given in equation \eqref{equation_RelativeEstErr}.}
\label{figure_RelativeEstErr_59x55}
\end{figure}

The reference value $\mathbf{x}_{\text{ref}}$ for measuring estimation errors is the true solution $\mathbf{x}_{\text{true}}$ or the LMMSE solution $\mathbf{x}_{\text{KF}}$ given by KF depending on the context. The true image as well as the reconstructed image of synthetic CO$_2$ plume are shown in Figure~\ref{figure_EnKFest_59x55} at low resolution ($59 \times 55$). The HiKF solution is also shown in Figure~\ref{figure_EnKFest_234x217} at high resolution ($234 \times 217$). The estimation errors relative to the true image as well as to the LMMSE estimate are plotted as a function of time in Figure~\ref{figure_RelativeEstErr_59x55}. As shown from the figures, HiKF reproduces the LMMSE estimates given by KF very accurately while EnKF produces noisy estimates even with an ensemble of $600$ realizations. The predictions become increasingly accurate over time as more data have been assimilated. Notice that after the CO$_2$ breakthrough (2 days after injection), the volume of CO$_2$ between the wells ceases to change. The crosswell seismic survey becomes less informative, which yields less improvement in estimation errors. 
 
\FloatBarrier

\begin{figure}
\centering\includegraphics{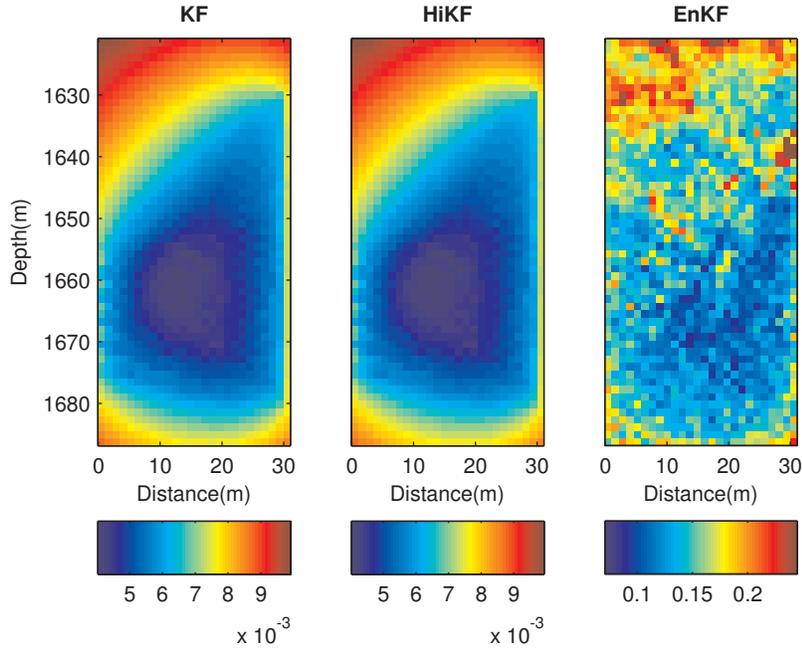}
\caption{Gridblock standard deviation of the compressional-wave slowness ($\times 10^{-4}s/m$) estimated by KF, HiKF and EnKF at last assimilation time step corresponding to Figure~\ref{figure_EnKFest_59x55}.}
\label{figure_UQ_final_result}
\end{figure}

Figure~\ref{figure_UQ_final_result} shows the estimated variance plotted as a function of the location at Day 5 after CO$_2$ injection. The uncertainty quantification given by KF indicates that the estimates inside the trapezoidal zone between two wells has less uncertainty compared to the remainder of the model domain, which coincides with the regions of high seismic ray coverage. As shown in the middle of the figure, HiKF reproduces the variance of KF very accurately. In contrast, the variance predicted by EnKF is $1000$ times greater than that given by KF and has to be plotted on a different color scale. As the regularization parameters and initial conditions are kept the same for the three experiments, the large variability of EnKF solutions can only result from sampling errors.    

\FloatBarrier

\begin{figure}
\centering\includegraphics{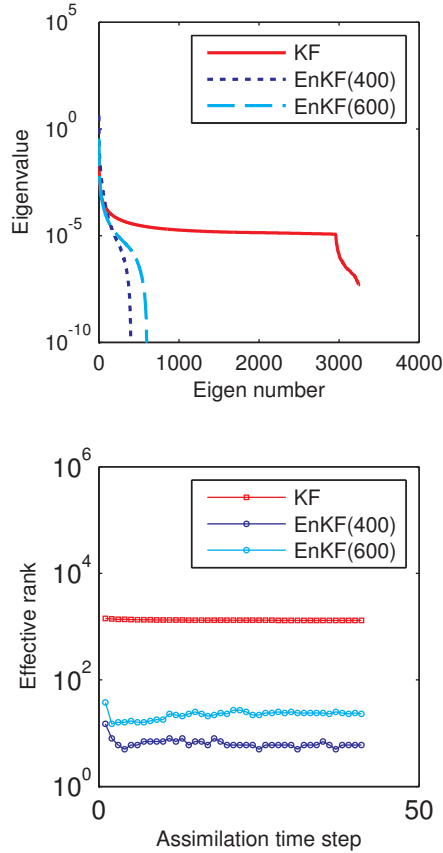}
\caption{Eigenvalue spectrums of the posterior covariance at last time step (top). Covariance given by EnKF are computed from 400- and 600-member ensembles respectively. Effective rank of the posterior covariance is plotted for each assimilation step (bottom).}
\label{figure_rank}
\end{figure}

We further investigate why EnKF cannot give a reliable uncertainty quantification in our case by plotting the eigenvalue spectrum and the effective rank of the posterior covariance in Figure~\ref{figure_rank}. In comparison with the eigenvalue spectrum at the last time step given by KF, which is flat and has a long tail, EnKF yields a steep eigenvalue spectrum with excess variance associated with the leading eigenvectors and insufficient variance with the rest. The same observation is reported in~\cite{hamill2001distance}. By increasing the ensemble size from $400$ to $600$, the variance in the tails of the spectrum increases yet the magnitude is still much lower than KF. Figure~\ref{figure_rank} also shows the effective rank of the posterior covariance for each assimilation time step. The effective rank is defined as the least number of principle eigenvalues needed to explain $95\%$ of the total variance. The plot suggests that EnKF suffers from rank deficiency associated with using insufficient ensemble size. Overall, the results show that the information embedded in the covariance of KF cannot be adequately represented by the proposed ensemble size. It is expected that with larger sample size, the variance given by EnKF will approach the values given by KF; however, this comes at the expense of increasing computational costs.

\FloatBarrier

\begin{figure}
\centering\includegraphics{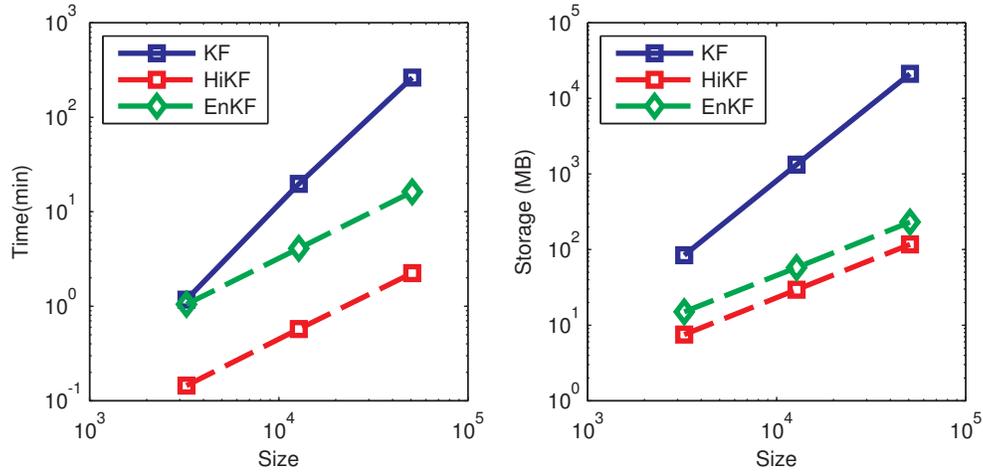}
\caption{Computational time (minutes) and storage cost (MB) of each data assimilation method plotted as a function of the number of unknowns. }
\label{figure_Cost}
\end{figure}

The computational cost of performing a single measurement update is summarized in Table~\ref{table_computationtimeC++} for each data assimilation method and plotted against the number of unknowns in Figure~\ref{figure_Cost} for comparison. All computations are performed using C++ on a PC with 2.40 GHz single-core CPU. Compared to the conventional KF, HiKF reduces the time of processing $~10^5$ unknowns from 4 hours to around 2 minutes, making it a promising tool for fast online data processing. As shown from Figure~\ref{figure_Cost}, the computational time and storage cost of KF grows quadratically with the number of unknowns, which involves propagating and storing a full error covariance matrix of size $m \times m$. In comparison, the overall cost of HiKF is $\mathcal{O}(m)$ because it operates on the cross-covariance matrix of size $m \times n$. The computational cost of HiKF comprises offline and online parts. The offline expense measures the cost of forming $QH^T$ from a dense covariance matrix $Q$, which needs to be computed only once in each monitoring event before any data is processed. The cost is normally $\mathcal{O}(m^2)$, but it can be reduced to $\mathcal{O}(m)$ by efficiently representing $Q$ as a $\mathcal{H}^2$ matrix~\cite{ambikasaran2013large}. The online portion only consists of forming the matrix-vector product associated with the cross-covariance, thus both the update and storage costs are $\mathcal{O}(m)$. The computational and storage costs of EnKF also scale linearly with the number of unknowns, as it propagates errors using an ensemble consisting of $N$ realizations of state vectors of size $m \times 1$ instead of a large covariance matrix. However, EnKF is computationally more expensive than HiKF in our case, as a very large ensemble size is required to produce accurate solutions.   

\begin{center}
\begin{table}
\center
\caption{Comparison of the HiKF with other KF variants}
\small{
\rowcolors{1}{}{gray!40}
\begin{tabular}{c|c|c|c}
\hline
\bf{Criteria} & \bf{KF} & \bf{HiKF} & \bf{EnKF}\\
\hline
Uncertainty propagation&Covariance & Cross-covariance & Ensembles \\
Estimation error & LMMSE & LMMSE & Large error for small ensemble size\\
Update cost & $\mathcal{O}(m^2)$ & $\mathcal{O}(m)$ & $\mathcal{O}(m)$\\
Storage cost  & $\mathcal{O}(m^2)$ & $\mathcal{O}(m)$ & $\mathcal{O}(m)$\\
Reliable risk analysis & Yes & Variance and cross-covariance & Require a large ensemble size\\
Generality & Linear forecast model & Random walk model & Any forecast model 
\end{tabular}
}
\label{table_comparison}
\end{table}
\end{center}

\begin{center}
\begin{table}
\center
\caption{Comparison of time taken on a 2.40 GHz single-core CPU for HiKF with other KF variants}
\small{
\rowcolors{1}{}{gray!40}
\begin{tabular}{c|c|c|c|c|c|c|c}
\hline
\bf{Grid size} & \multicolumn{1}{c}{\textbf{Precomputation (minutes)}} \vline & \multicolumn{3}{c}{\bf{Time for $41$ assimilation time steps (minutes)}} \vline & \multicolumn{3}{c}{\bf{Storage Cost (MB)}}\\
\hline
& HiKF & KF & EnKF & HiKF & KF & EnKF & HiKF \\
\hline  
$59 \times 55$ & $0.52$ &  $1.18$ &$1.05$&$0.14$ & 84.2 & 15.0 & 7.5  \\
$117 \times 109$ & $1.53$ & $19.8$ &$4.10$&$0.57$ & 1331.2 & 57.9 & 29.4  \\
$234 \times 217$ & $4.67$ & $263.7$ &$16.3$&$2.25$ & 21125.1 & 230.5 & 117  \\
Complexity&$\mathcal{O}(m)$&$\mathcal{O}(m^2)$&$\mathcal{O}(m)$&$\mathcal{O}(m)$& $\mathcal{O}(m^2)$& $\mathcal{O}(m)$ & $\mathcal{O}(m)$ \\
\end{tabular}
}
\label{table_computationtimeC++}
\end{table}
\end{center}

\begin{center}
\begin{table}
\center
\caption{Comparison of the HiKF with other KF variants}
\small{
\rowcolors{1}{}{gray!40}
\begin{tabular}{c|c|c|c}
\hline
\bf{Criteria} & \bf{KF} & \bf{HiKF} & \bf{EnKF}\\
\hline
Uncertainty propagation&Covariance & Cross-covariance & Ensembles \\
Estimation error & LMMSE & LMMSE & Large error for small ensemble size\\
Update cost & $\mathcal{O}(m^2)$ & $\mathcal{O}(m)$ & $\mathcal{O}(m)$\\
Storage cost  & $\mathcal{O}(m^2)$ & $\mathcal{O}(m)$ & $\mathcal{O}(m)$\\
Reliable risk analysis & Yes & Variance and cross-covariance & Require a large ensemble size\\
Generality & Linear forecast model & Random walk model & Any forecast model 
\end{tabular}
}
\label{table_comparison_1}
\end{table}
\end{center}
\FloatBarrier

\section{Conclusions}

We have presented a novel fast and accurate Kalman filter variant, HiKF, for solving quasi-continuous data assimilation problems. As shown in the numerical experiments, when vast data sets are collected in rapid succession, it is beneficial to adopt the random-walk forecast model and assume spatially-correlated Gaussian model errors described by a generalized covariance function. This allows the algorithm to explore the low-rank structure of the sub-blocks in the model error covariance matrix using $\mathcal{H}^2$-matrix algebra, and propagate errors without expensive operations involving a covariance matrix. As a result, the computational and storage costs of HiKF scale linearly with the number of unknowns.

A comprehensive comparison of HiKF, KF and EnKF can be found in Table~\ref{table_comparison}. As shown in the example case, by selecting a reasonable number of Chebyshev interpolation points, the solution given by HiKF can be as accurate as the LMMSE estimates given by KF within a controllable tolerance and with much less computational effort. To achieve the same accuracy, EnKF needs a high number of realizations in the example case. This is because when the ensemble size is small, the Monte Carlo approximation to the cross-covariance is low rank, which yields spurious covariance estimates and results in unphysical updates of the state variables. In comparison, HiKF propagates the cross-covariance in full rank, therefore gives accurate state and variance predictions.

The scope of this paper is to design a fast data-driven algorithm for quasi-continuous subsurface monitoring that is characterized by fast data acquisition using a permanent observation network. In such cases the random walk forecast model has certain advantages over a full-order forward model because it is fast, easy to implement, and the resulting solution relies more on the data and less on the model assumptions. Furthermore, we have shown in this paper that if the dynamical system can be described using the random walk state space model, then the Kalman filter solution can be obtained in linear running time using HiKF. However, for cases where data is sparse and can only be collected occasionally, the EnKF may be a better choice than the method presented here because the EnKF can incorporate a non-linear PDE solver for model forecasting. This allows one to rely less on data and more on model assumptions for correction. As data acquisition in the subsurface has became increasingly more affordable, collecting data in a quasi-continuous manner is both feasible and desirable. Thus, Kalman filtering can be used to continuously improve the quality of monitoring effects.  

\section{Acknowledgments}
This material is based upon work supported by ``US Department of Energy, National Energy Technology Laboratory" (DOE, NETL) under the award number \textendash DE-FE$0009260$: ``An Advanced Joint Inversion System for CO$_2$ Storage Modeling with Large Data Sets for Characterization and Real-Time Monitoring", and also by the ``National Science Foundation" \textemdash Division of Mathematical Sciences \textemdash under the award number: $1228275$. The author would like to thank Dr. Jonathan B. Ajo-Franklin, Thomas M. Daley, and Christine Doughty from the Lawrence Berkeley National Lab for sharing TOUGH2 and rock physics simulation data, and the support from ``The Global Climate and Energy Project" (GCEP) at Stanford.

\bibliographystyle{abbrv}
\bibliography{Kalman.bbl}

\begin{thebibliography}{10}

\bibitem{ajo2009optimal}
J.~B. Ajo-Franklin.
\newblock Optimal experiment design for time-lapse traveltime tomography.
\newblock {\em Geophysics}, 74(4):Q27--Q40, 2009.

\bibitem{ambikasaran2013fast}
S.~Ambikasaran and E.~Darve.
\newblock An $\mathcal{O}({N} \log {N})$ fast direct solver for partial
  hierarchically semi-separable matrices.
\newblock {\em Journal of Scientific Computing}, pages 1--25, 2013.

\bibitem{ambikasaran2013large}
S.~Ambikasaran, J.~Y. Li, P.~K. Kitanidis, and E.~Darve.
\newblock Large-scale stochastic linear inversion using hierarchical matrices.
\newblock {\em Computational Geosciences}, 17(6):913--927, 2013.

\bibitem{ambikasaran2013fast2}
S.~Ambikasaran, A.~K. Saibaba, E.~Darve, and P.~K. Kitanidis.
\newblock {\em Fast algorithms for Bayesian Inversion, The IMA Volumes in
  Mathematics and its Applications, vol. $156$}.
\newblock Springer-Verlag, New York, 2013.

\bibitem{anderson1979optimal}
B.~Anderson and J.~Moore.
\newblock {\em Optimal filtering}, volume~11.
\newblock Prentice-hall Englewood Cliffs, NJ, 1979.

\bibitem{anderson2001ensemble}
J.~L. Anderson.
\newblock An ensemble adjustment {K}alman filter for data assimilation.
\newblock {\em Monthly Weather Review}, 129(12):2884--2903, 2001.

\bibitem{anderson2007exploring}
J.~L. Anderson.
\newblock Exploring the need for localization in ensemble data assimilation
  using a hierarchical ensemble filter.
\newblock {\em Physica D: Nonlinear Phenomena}, 230(1):99--111, 2007.

\bibitem{arogunmati2009approach}
A.~Arogunmati and J.~M. Harris.
\newblock An approach for quasi-continuous time-lapse seismic monitoring with
  sparse data.
\newblock In {\em 79th Annual Meeting and International Exposition}, 2009.

\bibitem{bebendorf2000approximation}
M.~Bebendorf.
\newblock Approximation of boundary element matrices.
\newblock {\em Numerische Mathematik}, 86(4):565--589, 2000.

\bibitem{bebendorf2008hierarchical}
M.~Bebendorf.
\newblock {\em Hierarchical matrices: a means to efficiently solve elliptic
  boundary value problems}, volume~63.
\newblock Springer, 2008.

\bibitem{benson2006monitoring}
S.~Benson.
\newblock Monitoring carbon dioxide sequestration in deep geological formations
  for inventory verification and carbon credits.
\newblock In {\em SPE Annual Technical Conference and Exhibition}, 2006.

\bibitem{bishop2001adaptive}
C.~H. Bishop, B.~J. Etherton, and S.~J. Majumdar.
\newblock Adaptive sampling with the ensemble transform {Kalman} filter. {Part
  I: Theoretical} aspects.
\newblock {\em Monthly weather review}, 129(3):420--436, 2001.

\bibitem{borm2003hierarchical}
S.~B{\"o}rm, L.~Grasedyck, and W.~Hackbusch.
\newblock Hierarchical matrices.
\newblock {\em Lecture notes 21}, 2003.

\bibitem{borm2003introduction}
S.~B{\"o}rm, L.~Grasedyck, and W.~Hackbusch.
\newblock Introduction to hierarchical matrices with applications.
\newblock {\em Engineering Analysis with Boundary Elements}, 27(5):405--422,
  2003.

\bibitem{burgers1998analysis}
G.~Burgers, P.~van Leeuwen, and G.~Evensen.
\newblock Analysis scheme in the ensemble {K}alman filter.
\newblock {\em Monthly Weather Review}, 126:1719--1724, 1998.

\bibitem{chandrasekaran2006fast1}
S.~Chandrasekaran, P.~Dewilde, M.~Gu, W.~Lyons, and T.~Pals.
\newblock A fast solver for {HSS} representations via sparse matrices.
\newblock {\em SIAM Journal on Matrix Analysis and Applications}, 29(1):67--81,
  2006.

\bibitem{chandrasekaran2006fast}
S.~Chandrasekaran, M.~Gu, and T.~Pals.
\newblock A fast {ULV} decomposition solver for hierarchically semiseparable
  representations.
\newblock {\em SIAM Journal on Matrix Analysis and Applications},
  28(3):603--622, 2006.

\bibitem{daily1992electrical}
W.~Daily, A.~Ramirez, D.~LaBrecque, and J.~Nitao.
\newblock Electrical resistivity tomography of vadose water movement.
\newblock {\em Water Resources Research}, 28(5):1429--1442, 1992.

\bibitem{daley2007continuous}
T.~Daley, R.~Solbau, J.~Ajo-Franklin, and S.~Benson.
\newblock Continuous active-source seismic monitoring of injection in a brine
  aquifer.
\newblock {\em Geophysics}, 72(5):A57--A61, 2007.

\bibitem{evensen1994sequential}
G.~Evensen.
\newblock Sequential data assimilation with a nonlinear quasi-geostrophic model
  using monte carlo methods to forecast error statistics.
\newblock {\em Journal of Geophysical Research}, 99:10143--10162, 1994.

\bibitem{evensen2004sampling}
G.~Evensen.
\newblock Sampling strategies and square root analysis schemes for the enkf.
\newblock {\em Ocean Dynamics}, 54(6):539--560, 2004.

\bibitem{fahimuddin2010ensemble}
A.~Fahimuddin, S.~Aanonsen, and J.-A. Skjervheim.
\newblock Ensemble based 4d seismic history matching: Integration of different
  levels and types of seismic data.
\newblock In {\em SPE EUROPEC/EAGE Annual Conference and Exhibition}, 2010.

\bibitem{fong2009black}
W.~Fong and E.~Darve.
\newblock The black-box fast multipole method.
\newblock {\em Journal of Computational Physics}, 228(23):8712--8725, 2009.

\bibitem{Fritz:2009jn}
J.~Fritz, I.~Neuweiler, and W.~Nowak.
\newblock Application of {FFT}-based algorithms for large-scale universal
  kriging problems.
\newblock {\em Mathematical Geosciences}, 41(5):509--533, Apr. 2009.

\bibitem{fukumori1994approximate}
I.~Fukumori and P.~Malanotte-Rizzoli.
\newblock An approximate {K}alman filter for ocean data assimilation; an
  example with an idealized gulf stream model.
\newblock {\em Journal of Geophysical Research}, 1994.

\bibitem{goreinov1997theory}
S.~A. Goreinov, E.~E. Tyrtyshnikov, and N.~L. Zamarashkin.
\newblock A theory of pseudoskeleton approximations.
\newblock {\em Linear Algebra and its Applications}, 261(1):1--21, 1997.

\bibitem{greengard1987fast}
L.~Greengard and V.~Rokhlin.
\newblock A fast algorithm for particle simulations.
\newblock {\em Journal of computational physics}, 73(2):325--348, 1987.

\bibitem{guivant2001optimization}
J.~E. Guivant and E.~M. Nebot.
\newblock Optimization of the simultaneous localization and map-building
  algorithm for real-time implementation.
\newblock {\em Robotics and Automation, IEEE Transactions on}, 17(3):242--257,
  2001.

\bibitem{hackbusch1999sparse}
W.~Hackbusch.
\newblock A sparse matrix arithmetic based on {$\mathcal{H}$}-matrices. part
  {I}: Introduction to {$\mathcal{H}$}-matrices.
\newblock {\em Computing}, 62(2):89--108, 1999.

\bibitem{hackbusch2000sparse}
W.~Hackbusch and B.~Khoromskij.
\newblock A sparse {$\mathcal{H}$}-matrix arithmetic. part {II}: Application to
  multi-dimensional problems.
\newblock {\em Computing}, 64(1):21--47, 2000.

\bibitem{hamill2001distance}
T.~Hamill, J.~Whitaker, and C.~Snyder.
\newblock Distance-dependent filtering of background error covariance estimates
  in an ensemble {K}alman filter.
\newblock {\em Monthly Weather Review}, 129(11):2776--2790, 2001.

\bibitem{houtekamer1998data}
P.~Houtekamer and H.~Mitchell.
\newblock Data assimilation using an ensemble {K}alman filter technique.
\newblock {\em Monthly Weather Review}, 126(3):796--811, 1998.

\bibitem{houtekamer2001a}
P.~Houtekamer and H.~Mitchell.
\newblock A sequential ensemble {K}alman filter for atmospheric data
  assimilation.
\newblock {\em Monthly Weather Review}, 129:123--137, 2001.

\bibitem{houtekamer2001sequential}
P.~L. Houtekamer and H.~L. Mitchell.
\newblock A sequential ensemble {K}alman filter for atmospheric data
  assimilation.
\newblock {\em Monthly Weather Review}, 129(1):123--137, 2001.

\bibitem{hubbard2001hydrogeological}
S.~Hubbard, J.~Chen, J.~Peterson, E.~Majer, K.~Williams, D.~Swift, B.~Mailloux,
  and Y.~Rubin.
\newblock Hydrogeological characterization of the south oyster bacterial
  transport site using geophysical data.
\newblock {\em Water Resources Research}, 37(10):2431--2456, 2001.

\bibitem{Kalman:1960tn}
R.~Kalman.
\newblock A new approach to linear filtering and prediction problems.
\newblock {\em Journal of Basic Engineering}, 82(1):35--45, 1960.

\bibitem{kepert2004ensemble}
J.~D. Kepert.
\newblock On ensemble representation of the observation-error covariance in the
  ensemble {K}alman filter.
\newblock {\em Ocean Dynamics}, 54(6):561--569, 2004.

\bibitem{kim20094}
J.~Kim, M.~Yi, S.~Park, and J.~Kim.
\newblock {4-D} inversion of {DC} resistivity monitoring data acquired over a
  dynamically changing earth model.
\newblock {\em Journal of Applied Geophysics}, 68(4):522--532, 2009.

\bibitem{pk1997introduction}
P.~Kitanidis.
\newblock {\em Introduction to geostatistics: applications in hydrogeology}.
\newblock Cambridge University Press, 1997.

\bibitem{kitanidis1983statistical}
P.~K. Kitanidis.
\newblock Statistical estimation of polynomial generalized covariance functions
  and hydrologic applications.
\newblock {\em Water Resources Research}, 19(4):909--921, 1983.

\bibitem{kitanidis1993generalized}
P.~K. Kitanidis.
\newblock Generalized covariance functions in estimation.
\newblock {\em Mathematical Geology}, 25(5):525--540, 1993.

\bibitem{kitanidis1999generalized}
P.~K. Kitanidis.
\newblock Generalized covariance functions associated with the laplace equation
  and their use in interpolation and inverse problems.
\newblock {\em Water Resources Research}, 35(5):1361--1367, 1999.

\bibitem{lazaratos1996crosswell}
S.~Lazaratos and B.~Marion.
\newblock Crosswell seismic imaging of reservoir changes caused by {CO}$_2$
  injection.
\newblock {\em The Leading Edge}, 16(9):1300--1308, 1997.

\bibitem{martinsson2005fast}
P.-G. Martinsson and V.~Rokhlin.
\newblock A fast direct solver for boundary integral equations in two
  dimensions.
\newblock {\em Journal of Computational Physics}, 205(1):1--23, 2005.

\bibitem{nenna2011application}
V.~Nenna, A.~Pidlisecky, and R.~Knight.
\newblock Application of an extended {K}alman filter approach to inversion of
  time-lapse electrical resistivity imaging data for monitoring recharge.
\newblock {\em Water Resources Research}, 47(10), 2011.

\bibitem{neumaier1998solving}
A.~Neumaier.
\newblock Solving ill-conditioned and singular linear systems: A tutorial on
  regularization.
\newblock {\em Siam Review}, 40(3):636--666, 1998.

\bibitem{nowak2003efficient}
W.~Nowak, S.~Tenkleve, and O.~Cirpka.
\newblock Efficient computation of linearized cross-covariance and
  auto-covariance matrices of interdependent quantities.
\newblock {\em Mathematical Geology}, 35(1):53--66, 2003.

\bibitem{Pnevmatikakis:2013jh}
E.~A. Pnevmatikakis, K.~R. Rad, J.~Huggins, and L.~Paninski.
\newblock {Fast Kalman filtering and forward-backward smoothing via a low-rank
  perturbative approach}.
\newblock {\em Journal of Computational and Graphical Statistics}, Jan. 2013.

\bibitem{pruess1991tough2}
K.~Pruess.
\newblock {TOUGH}2: A general-purpose numerical simulator for multiphase fluid
  and heat flow.
\newblock {\em NASA STI/Recon Technical Report N}, 92(14316), 1991.

\bibitem{quan2008stochastic}
Y.~Quan and J.~Harris.
\newblock Stochastic seismic inversion using both waveform and traveltime data
  and its application to time-lapse monitoring.
\newblock In {\em 2008 SEG Annual Meeting}, 2008.

\bibitem{rjasanow2002adaptive}
S.~Rjasanow.
\newblock Adaptive cross approximation of dense matrices.
\newblock In {\em Int. Association Boundary Element Methods Conf., IABEM},
  pages 28--30, 2002.

\bibitem{saibaba2012application}
A.~Saibaba, S.~Ambikasaran, J.~Li, P.~Kitanidis, and E.~Darve.
\newblock Application of hierarchical matrices to linear inverse problems in
  geostatistics.
\newblock {\em OGST Revue d'IFP Energies Nouvelles}, 67(5):857--875, 2012.

\bibitem{starks1982estimation}
T.~Starks and J.~Fang.
\newblock On the estimation of the generalized covariance function.
\newblock {\em Journal of the International Association for Mathematical
  Geology}, 14(1):57--64, 1982.

\bibitem{Tippett:2003vi}
M.~K. Tippett, J.~L. Anderson, C.~H. Bishop, T.~M. Hamill, and J.~S. Whitaker.
\newblock Ensemble square root filters*.
\newblock {\em Monthly Weather Review}, 131(7):1485--1490, 2003.

\bibitem{TuanPham:1998hg}
D.~Tuan~Pham, J.~Verron, and M.~Christine~Roubaud.
\newblock {A singular evolutive extended {K}alman filter for data assimilation
  in oceanography}.
\newblock {\em Journal of Marine Systems}, 16(3-4):323--340, Oct. 1998.

\bibitem{vauhkonen1998kalman}
M.~Vauhkonen, P.~Karjalainen, and J.~Kaipio.
\newblock A {Kalman} filter approach to track fast impedance changes in
  electrical impedance tomography.
\newblock {\em Biomedical Engineering, IEEE Transactions on}, 45(4):486--493,
  1998.

\bibitem{wang1998seismic}
Z.~Wang, M.~E. Cates, and R.~T. Langan.
\newblock Seismic monitoring of a {CO}$_2$ flood in a carbonate reservoir: A
  rock physics study.
\newblock {\em Geophysics}, 63(5):1604--1617, 1998.

\bibitem{White:2001ul}
J.~White.
\newblock Computed seismic speeds and attenuation in rocks with partial gas
  saturation.
\newblock {\em Geophysics}, 40(2):224--232, Sept. 2001.

\end{thebibliography}

\end{document}